# Optimal Design and Operation of Integrated Microgrids under Intermittent Renewable Energy Sources Coupled with Green Hydrogen and Demand Scenarios

Su Meyra Tatar[a], Handan Akulker[a], Hasan Sildir[b] and Erdal Aydin[a,c,d*]

[a]Department of Chemical Engineering, Bogazici University, Bebek, Istanbul 34342, Turkey
[b]Department of Chemical Engineering, Gebze Technical University, Kocaeli 41400, Turkey
[c]Department of Chemical and Biological Engineering, Koç University, Istanbul 34457, Turkey
[d]Koç University TUPRAS Energy Center (KUTEM), Koç University, Istanbul, 34450, Turkey

*eaydin@ku.edu.tr

**Abstract**

Investigation of energy systems integrated with green chemical conversion, and in particular combination of green hydrogen and synthetic methanation, is still a scarce subject in the literature in terms of optimal design and operation for energy grids under weather intermittency and demand uncertainty. In this work, a multi-period mixed-integer linear programming (MILP) model is formulated to identify the optimal design and operation of integrated energy grids including such chemical conversion systems. Under current carbon dioxide limitations, this model computes the best configuration of the renewable and non-renewable-based generators, from a large candidate pool containing thirty-nine different equipment, their optimal rated powers, capacities and scheduling sequences. Three different scenarios are generated for a specific location. We observed that photovoltaic, oil co-generator, reciprocating ICE, micro turbine, and bio-gasifier are the equipment that is commonly chosen under the three different scenarios. Results also show that concepts such as green hydrogen and power-to-gas are currently not preferable for the investigated location.

# 1. Introduction

Energy producers are globally struggling with the inadequacy of energy resources and energy efficiency problems, which have become even more active in recent years due to COVID-19 pandemics. For example, the ongoing natural gas crisis, which has emerged in Turkey and most of the other countries in late 2021, looks sure to continue in the coming future (Gilbert et al., 2021). Hence, it is crucial to find new ways to produce energy sources, such as natural gas, hydrogen, etc. *Power to Gas* or *Power to X* and *Green Hydrogen* technologies have been developed for this purpose. However, they are required to be integrated with existing energy hub technologies to secure the supply. Moreover, their operational plans should be simultaneously carried out with the hubs' energy production plan.

Another issue is energy efficiency, a significant economic growth and prosperity metric. In particular, optimal energy production, distribution, and consumption have vital importance on the economic independence. For instance, Turkey's account deficit is primarily caused by $60 billion in energy imports annually (Orhan and Nergiz, 2014). On the other hand, Turkey is one of the richest countries in terms of renewable energy resources. At the same time, Turkey has a rapidly expanding population. Thus, it is essential to integrate renewable energy-based equipment into energy hubs together with *Power to Gas* and *Green Hydrogen Technologies*, considering current and future trends.

In addition to those, energy producers have been struggling with the effects of greenhouse gases. In the countries which have signed and ratified the Paris Agreement, energy producers are obliged to fulfill its sanctions (The Paris Agreement, 2015). Turkey is one of these countries, and it has no certain sanctions for energy producers to reduce carbon dioxide emission, such as the cap-and-trade system or emission taxing, yet.

In the literature, there are several studies on MILP-based integrated energy hub optimization handling a few of the aforementioned issues. In (Eladl et al., 2020), optimum operation and configuration of a hybrid renewable energy integration system are studied to minimize $CO_2$ emission by a Genetic Algorithm, without using rigorous and derivative-based solution methods. In (Nicolosi et al., 2021), the optimal unit commitment of a micro-grid is performed by MILP to reduce the impact of greenhouse gas emissions, without taking into account the optimal design phase simultaneously. In (Naderi et al., 2022), Demand Response based MILP is formulated for optimal energy management and sizing of a hybrid microgrid under wind intermittency only. Green hydrogen technology is utilized to obtain hydrogen by water electrolysis with excess electricity accumulated in an energy hub. In (Tostado-Véliz et al., 2022), a novel optimal scheduling model for robust optimal scheduling of energy hubs, including green hydrogen technology, is formulated by MILP. This study also evaluates the impact of green hydrogen technology on the reduction of fossil fuels consumption. In (Nasiri et al., 2021), a bi-level optimization is conducted to assess the strategic behavior of a multi-energy system consisting of green hydrogen technology. In (Silvente et al., 2017a), the optimal design of off-grid hydrogen-battery energy systems is done by MILP. Power to gas is a relatively new technology that connects power and gas networks by converting power to natural gas by electrolysis and methanation, respectively.

In (Ikäheimo et al., 2022), optimal scheduling of energy hubs considering a power-to-gas system is studied using MILP. The financial effects of power to gas systems are evaluated. In (Yuan et al., 2020), probabilistic scheduling of a hybrid energy system including power-to-gas technology is performed by stochastic MILP handling the uncertainties related to the electrical loads, wind power, and electricity price. (Dong et al., 2022) proposes an optimal energy production plan for an integrated energy hub including power to gas technology and gas-fired units with carbon capture systems using MILP. It is again vital to state that aforementioned studies do not cover the optimal design phase. On the other hand, in order to fully assess the impact of such current trends, design phase must be taken into account.

In this work, a multi-period mixed-integer linear programming (MILP) model is formulated to identify the optimal design and operation of integrated energy grids under weather intermittency and demand uncertainty considering the combination of green hydrogen and synthetic methanation as a candidate. This multi-period model selects the equipment to minimize the net present cost (total project cost over the lifetime of the microgrid) of the integrated microgrid for optimal design and operational plan by considering the energy issues mentioned above. For example, in the suggested decision-making model, carbon dioxide emission is limited according to the Paris Agreement. In addition, *Power to Gas* and *Green Hydrogen* technologies are included in the candidate equipment pool. Formulated MILP model combines the optimal design and operation on the same level in a simultaneous fashion; it computes both optimal investment decisions in terms of equipment selection and sizing and the optimal sub-hourly scheduling sequences of this equipment. This way, optimal design and operation are investigated simultaneously under weather intermittency and demand uncertainty combined with chemical conversion such as green hydrogen production and synthetic methanation. Moreover, linking and logical constraints between the binary and continuous variables tighten the search space for the MILP solvers and therefore reduce the overall computational load of the model. In addition, this model includes features as disaggregation of bilinear terms to linearize the generic nonlinear unit commitment equations in an exact manner once these equipment are decided to be installed (Yokoyama and Ito, 2004). Such improvements result in the ability to formulate the problem as a mixed-integer linear programming (MILP), which is relatively easier to solve using solvers CPLEX or Gurobi, instead of derivative-free heuristic algorithms. Finally, this multi-period MILP model considers peak shaving and spinning reserves, which are also crucial metrics for energy management and should be taken into account in the design phase rather than the scheduling or control (Takriti et al., 1996).

Additionally, this work also compares the results of three different case studies based on statistical scenarios, including uncertainties coming from renewable energy sources and electricity demand. To generate related and dependable scenarios, a practical scenario generation method using actual data is proposed. This scenario generation method consists of three main consecutive procedures. The first one is the formulation of empirical mathematical models for wind speed, air temperature, solar irradiance, population, and electricity consumption, which are trained using historical data. Only a few scenario generation methods, used for energy optimization in general, provide the flexibility to formulate empirical models with parameters for uncertainty sources (Ehsan and Yang, 2019; Farrokhifar et al., 2020; Osório et al., 2015; Tanaka and Ohmori, 2016; Yu et al., 2019). One exception for empirical models would be the proposition of the Weibull distribution for the wind speed, which is a worldwide generalization. The second procedure of the proposed scenario generation method is to propagate the parameter uncertainty by using the covariance matrix and Taylor series approximation. This widens the prediction intervals of probability distribution functions, which in turn improves the reliability and robustness of the generated scenarios. In the literature, Monte Carlo simulation is the common approach for the probabilistic uncertainty modeling. Compared with the Taylor series approximation, Monte Carlo simulation is computationally more expensive due to its iterative behavior (Carpinelli et al., 2015; Zakaria et al., 2020). Finally, the last step is inverse transform sampling. This procedure provides sampling from a predetermined interval under desired likelihood, such as likely, mid-likely, and unlikely, whose main advantage is that a scenario reduction method is not required, which would easily bring additional computational effort (Beltran-Royo, 2019; Silvente et al., 2017b; Yildiz et al., 2017).

Accordingly, the key contributions of this study can be summarized as:

1) An MILP model is formulated to identify the optimal design and operation of integrated energy grids under weather intermittency and demand uncertainty. A practical probabilistic scenario generation method is proposed to handle the uncertainties related to wind speed, solar irradiance, air temperature, population, and electricity demand of a city. This method also provides specific model equations for the uncertainty sources in a different way from most of the studies related to scenario generation in renewable energy grids (Ehsan and Yang, 2019; Farrokhifar et al., 2020; Osório et al., 2015; Tanaka and Ohmori, 2016; Yu et al., 2019).

2) The variety and specifications of candidate equipment are determined to follow the newest trends in the energy sector, such as green hydrogen and power-to-gas technology in terms of optimal design and operation simultaneously. In addition, a $CO_2$ limit is set according to Paris Agreement (The Paris Agreement, 2015).

3) The case studies including likely, mid-likely, and unlikely scenarios, provide a beneficial comparison to understand the effect of uncertainties coming from renewable sources and demand.

The paper is structured as follows: Section 2 details the proposed MILP decision-making model. Section 3 presents the suggested scenario generation algorithm. Section 4 shows the case study results in terms of optimal design and operation for three different scenarios, including the discussions. Finally, Section 5 concludes this study.

## 2. The MILP Decision Making Model

| Nomenclature | | | |
|---|---|---|---|
| Sets | | Parameters (Cont'd) | |
| $I$ | Set of equipment, $i$ | $LCO2$ | Daily $CO_2$ emission limit, $gCO_2/kW$ |
| $G$ | Set of generators ($G \in I$) | $rp_i^{min}$ | Minimum introducible rated power of equipment $i$, |
| $S$ | Set of storage units ($S \in I$) | $rp_i^{max}$ | Maximum introducible rated power of equipment $i$ |
| $R$ | Set of renewables ($R \in I$) | $b_i^{min}$ | Minimum introducible capacity of storage ($i \in S$) |
| $K$ | Set of years, $k$ | $b_i^{max}$ | Maximum introducible capacity of storage ($i \in S$) |
| $T$ | Set of time-intervals, $t$ | $p_{ikt}^{min}$ | Minimum operating power fraction of equipment $i$ |
| $N$ | Set of resources, $n$ | $p_{ikt}^{max}$ | Maximum operating power fraction of equipment $i$ |
| Decision Variables | | $q_{ikt}^{min}$ | Minimum operating charged energy fraction |
| $rp_i$ | Rated power of equipment $i$ | $q_{ikt}^{max}$ | Maximum operating charged energy fraction |
| $b_i$ | Capacity of $i$, ($i \in S$) | $g_{in}$ | Generation of resource $n$ |
| $a_i$ | Binary decision to install | $c_{in}$ | Consumption of resource $n$ |
| $sp_{nk}$ | Spinning reserve for resource $n$ | $d_{nkt}$ | Demand of resource $n$ |
| $kc_{ikt}$ | Binary decision to run a generator | $\bar{p}_{ikt}$ | Operating power coefficient |
| $ks_{ikt}$ | Binary decision to run a storage | $\alpha_i^0$ | Unit cost for building equipment $i$ |
| $p_{ikt}$ | Operating power of generator $i$ | $\beta_i^0$ | Unit cost for building a storage unit ($i \in S$) |
| $pch_{ikt}$ | Charging power of equipment | $\gamma_i^0$ | Fixed cost for building equipment $i$ |
| $pdch_{ikt}$ | Discharging power of equipment | $\alpha_i^k$ | Unit cost for maintenance of equipment $i$ |
| $soc_{ik0}$ | Initially stored energy | $\beta_i^k$ | Unit cost for maintenance of a storage unit ($i \in S$) |
| $soc_{ikt}$ | Reserved energy of a storage device | $\gamma_i^k$ | Fixed cost for maintenance of equipment $i$ in year $k$ |
| $yx_{nkt}$ | Surplus output of resource $n$ | $\phi_{nkt}^+$ | Cost for unit system output of resource $n$ |
| $u_{nkt}$ | System input of resource $n$ | $\phi_{nkt}^-$ | Cost for unit system input of resource $n$ |
| $\xi$ | Peak penalty | $\delta$ | Peak penalty coefficient |
| Parameters | | $\eta$ | Efficiency of the photovoltaic cell |
| $inf$ | Inflation rate | $\Phi$ | Solar irradiance, in kW/area |
| $D$ | 365 days per year | $\kappa$ | Temperature correction factor |
| $\Delta T$ | Time interval, in hour | $T_p$ | Ambient temperature, in $^0C$ |
| $N_I$ | Maximum introducible equipment | $T_{p,ref}$ | Reference temperature, in $^0C$ |
| $v_{cut,in}$ | Cut-in wind speed, in m/sec | $v_{rated}$ | Rated wind speed, in m/sec |
| $v_{cut,out}$ | Cut-out wind speed, in m/sec | $v$ | Actual wind speed, in m/sec |

## 2.1. Objective function

The objective function is to minimize the net present cost of the integrated microgrid over its lifetime, 20 years (Kharrich et al., 2021):

$$\min \quad f^{Initial}(rp, b, a) + \sum_{k \in \mathcal{K}} f_k^{O\&M}(rp, b, a, yx_{nkt}, u_{nkt}) \tag{2.1}$$

The fundamental assumption of the proposed model is the quasi-steady-state for each time step. Eqn. 2.1. represents the optimization cost function which is to minimize the overall cost.

## 2.2. Installation and sizing decisions

$$\sum_{i \in I} a_i \leq N_I, \tag{2.2}$$

$$rp_i^{min} a_i \leq rp_i \leq rp_i^{max} a_i \; ; \quad \forall \, i \in G, \tag{2.3}$$

$$b_i^{min} a_i \leq b_i \leq b_i^{max} a_i \; ; \quad \forall \, i \in S, \tag{2.4}$$

$$a_i \in \{0,1\} \; ; \forall \, i \in I.$$

$a$ is the binary variable representing the installation decision of any equipment of set $I$. Eqn. 2.2 stands for the fact that the total number of installable equipment is limited by $N_I$, which is determined as ten due to the land constraint. Eqn. 2.3. shows the rated power, $rp$, constraints of any generator. Similarly, Eqn. 2.4 represents the capacity constraints of any storage unit where $b$ is the capacity of a particular storage unit. Please note that these equations also behave as logical constraints since the sizing decisions, which can be rated power and storage capacity depending on the type of the equipment, are enforced to be zero if the corresponding installation decision is not made. Accordingly, this formulation reduces the solution effort for the MILP problem extensively.

## 2.3. Unit Commitment and Operational Constraints

$$p_{ikt}^{min}(rp_i - (1 - kc_{ikt})rp_i^{max}) \leq p_{ikt} \leq p_{ikt}^{max} rp_i \; ; \quad \forall i \in G, \forall k \in \mathcal{K}, \quad \forall t \in \mathcal{T} \tag{2.5}$$

$$0 \leq p_{ikt} \leq p_{ikt}^{max} rp_i^{max} kc_{ikt} \; ; \quad \forall i \in G, \quad \forall k \in \mathcal{K}, \quad \forall t \in \mathcal{T} \tag{2.6}$$

$$kc_{ikt} \in \{0,1\} \; ; \quad \forall i \in G, \quad \forall k \in \mathcal{K}, \quad \forall t \in \mathcal{T}$$

Eqns. 2.5 and 2.6 show the installed generators' operational upper and lower power limits when switched on and switched off. $kc_{ikt}$ corresponds to the binary decision to operate a generator or not. Hence, while $kc_{ikt}$ equals to '1' corresponds to switched on mode, $kc_{ikt}$ equals to '0' corresponds to switched off mode. Unless an arbitrary installation decision, $a_i$, is made, the rated power, $rp_i$, is set automatically to zero, making the operational power equal to zero. Conversely, if the installation decision is made, then the rated power and the operational power of the equipment are bounded by the lower and upper limits. Furthermore, please note that these capacity equations are obtained by disaggregating the bilinear terms ($rp_i kc_{ikt}$), which are continuous and binary, followed by the exact linearization. The same formulation is also applied to the battery limit terms. Detailed information about this trick can be found in (Yokoyama and Ito, 2004), where installation decisions are not present.

$$0 \leq pch_{ikt} \leq p_{ikt}^{max} rp_i ; \quad \forall i \in S, \quad \forall k \in \mathcal{K}, \quad \forall t \in \mathcal{T} \quad (2.7)$$

$$0 \leq pdch_{ikt} \leq p_{ikt}^{max} rp_i ; \quad \forall i \in S, \quad \forall k \in \mathcal{K}, \quad \forall t \in \mathcal{T} \quad (2.8)$$

$$pch_{ikt} \leq p_{ikt}^{max} rp_i^{max} ks_{ikt} ; \quad \forall i \in S, \quad \forall k \in \mathcal{K}, \quad \forall t \in \mathcal{T} \quad (2.9)$$

$$pdch_{ikt} \leq p_{ikt}^{max} rp_i^{max} (1 - ks_{ikt}) ; \quad \forall i \in S, \quad \forall k \in \mathcal{K}, \quad \forall t \in \mathcal{T} \quad (2.10)$$

$$ks_{ikt} \in \{0,1\} ; \quad \forall i \in S, \quad \forall k \in \mathcal{K}, \quad \forall t \in \mathcal{T}$$

Equations from 2.7 to 2.10 show the charging and discharging power limits of a storage device. $p_{ikt}^{max}$ is a dimensionless number that is the maximum operating power of equipment by the fraction of the rated power. Hence, Eqns. 2.7-2.8 show that the $pch_{ikt}$ and $pdch_{ikt}$ variables, namely the charging and discharging power of the storage device, are both bounded with the rated power of the equipment. Please note that $ks_{ikt}$ corresponds to the binary decisions to run a storage device. Thus, Eqns. 2.9 and 2.10 correspond to the constraint that the storage devices cannot be charged and discharged simultaneously.

## 2.4. Storage Constraints

$$q_{ikt}^{min} b_i \leq soc_{ikt} \leq q_{ikt}^{max} b_i ; \forall i \in S, \quad \forall k \in \mathcal{K}, \quad \forall t \in \mathcal{T} \quad (2.11)$$

$$soc_{ik|T-1|} = soc_{ik0} ; \forall i \in S, \quad \forall k \in \mathcal{K} \quad (2.12)$$

$$soc_{ikt} = \begin{cases} soc_{ik0} + \Delta T(pch_{ikt} - pdch_{ikt}), & if \ t = t(1) \\ soc_{ikt-1} + \Delta T(pch_{ikt} - pdch_{ikt}), & otherwise \end{cases} \quad (2.13)$$
$$\forall i \in S, \quad \forall k \in \mathcal{K}, \quad \forall t \in \mathcal{T}$$

Reserved energy in a storage unit at the end of the time interval is represented by $soc_{ikt}$, which should be in between 20% and 80% of the storage capacity due to safety and equipment life-related considerations, as stated in Eqn. 2.11. $q_{ikt}^{min}$ and $q_{ikt}^{max}$ are the minimum and maximum operating charged energy of a storage unit. Moreover, Eqn. 2.11-2.13 represent capacity limits of a storage device, the fact that initial and final conditions of the storage device in a particular day must be the same and charging and discharging of a storage device, respectively.

## 2.5. Non-negativity Constraints for System Inputs and Outputs

$$u_{nkt} \geq 0; \quad \forall n \in \mathcal{N}, \quad \forall k \in \mathcal{K}, \quad \forall t \in \mathcal{T} \quad (2.14)$$

$$yx_{nkt} \geq 0; \quad \forall n \in \mathcal{N}, \quad \forall k \in \mathcal{K}, \quad \forall t \in \mathcal{T} \quad (2.15)$$

While the $u_{nkt}$ term represents the system input of resources, the $yx_{nkt}$ term represents the surplus output of resource *n*.

## 2.6. Material and Power Balance

$$\sum_{i \in G} g_{in} p_{ikt} + \sum_{i \in S} g_{in} pdch_{ikt} + \sum_{i \in \mathcal{R}} g_{in} \bar{p}_{ikt} rp_i + u_{nkt} = \sum_{i \in G} c_{in} p_{ikt}$$
$$+ \sum_{i \in S} c_{in} pch_{ikt} + \sum_{i \in \mathcal{R}} c_{in} \bar{p}_{ikt} rp_i + sp_{nk} + yx_{nkt} + d_{nkt} \quad (2.16)$$
$$\forall n \in \mathcal{N}, \quad \forall k \in \mathcal{K}, \quad \forall t \in \mathcal{T}$$

Total generation, consumption, input, output, spin, and demand balances for the resources are shown by Eqn. 2.16. $g_{in}$ and $c_{in}$ represent the generation and consumption of any equipment *I* for any resource *n*. Corresponding values are given in Tables A-3 and A-4 in the appendix. Eqn. 2.16 equates the sum of generation and system input of resource *n* to the sum of consumption, surplus output, and demand of resource *n* and spinning reserve of a generator producing resource *n*. The spinning reserve must be less than three percent of the maximum demand throughout the day (Takriti et al., 1996).

Some of the non-renewable energy-based equipment are CHP (combined heat and power) units where both heat and electricity are produced. Additionally, some of the generators and renewable equipment produce only power. Ideally, generators can convert the fuel they consume into electricity or electricity and heat without any loss. The data in generation and consumption tables, given in Appendix Tables of A.3 and A.4, are obtained from the literature. They represent the generation or consumption of resource *n* when equipment *i* is run at the unit power (kW) for unit time (hour). The data in the tables are represented with $g_{in}$ and $c_{in}$ as variable in equation 2.16, while '*i*' represents the set of equipment and '*n*' represents the set of resources.

## 2.7. Cost Calculations

$$f^{initial}(rp, b, a) = \sum_{i \in I}(\alpha_i^0 rp_i + \gamma_i^0 a_i) + \sum_{i \in S} \beta_i^0 b_i \tag{2.17}$$

$$f_k^{O\&M} = f_k^{Operational} + f_k^{Maintenance} \tag{2.18}$$

$$f_k^{Operational} = D\Delta T \sum_{t \in T} \sum_{n \in \mathcal{N}} (1 + inf)^{k-1}(\phi_{nkt}^+ yx_{nkt} + \phi_{nkt}^- u_{nkt}) + D|\mathcal{K}||T|\xi \tag{2.19}$$

$$f_k^{Maintenance} = \sum_{i \in I}(1 + inf)^{k-1}(\alpha_i^k rp_i + \gamma_i^k a_i) + \sum_{i \in S}(1 + inf)^{k-1}\beta_i^k b_i \tag{2.20}$$

We should note that the income from the cap and trade of $CO_2$ is not considered in this study. Cap and trade prices might be fairly negligible for the case studies shown in this work. However, these prices are expected to increase in the future and fairly affect the optimal results for larger scale studies.

Eqn. 2.17 and 2.18 define the initial and O&M costs. O&M costs are the summation of total operational and maintenance costs in year *k* as shown in Eqn. 2.19 and 2.20. Cost terms are found with the multiplication of cost with corresponding rated power or capacity value, where inflation ($inf$) is also taken into account. $\phi_{nkt}^- u_{nkt}$ term represents the raw material cost where $\phi_{nkt}^-$ is the cost for unit system input of resource *n*. $\phi_{nkt}^+ yx_{nkt}$ term can be considered as penalty term for producing excess resource *n*. Although this part is neglected in the case studies, please also note that related cost terms can be easily added to the objective function using the proposed formulation. The last term in Eqn. 2.21 is the penalty cost term that is associated with peak shaving. Here, $|T|$ stands for the size of time intervals and *D* is the number of days in a year.

$$\xi \geq \delta(p_{gkt})(g_{,Electricity,kt}) \tag{2.21}$$

Eqn. 2.21 stands for the definition of the peak penalty selection equation where $\xi$ and $\delta$ stand for peak penalty and peak penalty coefficient, respectively. Please note that the peak penalty is computed via the decision-making model and must be larger than the maximum electricity production from the generators. This way, peak supply from the generators can be reduced or shaved. The main advantages of this formulation are that it is a linear formulation, and does not require a non-smooth function, e.g., max operator, as opposed to many other conventional approaches (Feng et al., 2019). We provide relatively smaller values for the peak penalty coefficients as peak shaving is not the main focus of this study. On the other hand, analyzing the values of price coefficients on the optimal results and solution times is an active research area.

## 2.8. Emission Constraints

$$\sum_{t \in T}(yx_{,CO2,kt}) \leq LCO2 \tag{2.22}$$

Eqn. 2.22 defines the daily $CO_2$ emission limit.

## 2.9. Renewable Equipment Generations

$$\bar{p}_{,PV,kt} = \eta \Phi \left(1 - \kappa(T_p - T_{p,ref})\right) \tag{2.23}$$

$$\bar{p}_{,WT,kt} = \begin{cases} 0 & if\ v < v_{cut,in} \\ 0 & if\ v > v_{cut,out} \\ 1 & if\ v_{rated} \leq v \leq v_{cut,out} \\ \dfrac{v - v_{cut,in}}{v_{rated} - v_{cut,in}} & if\ v_{cut,in} \leq v < v_{rated} \end{cases} \tag{2.24}$$

Eqn. 2.23 explains the operating power coefficient of a photovoltaic cell directly proportional to the solar irradiance, $\Phi$, and efficiency of the photovoltaic cell, $\eta$. The operating power of the wind turbines depends on the actual wind speed, $v$. If the wind velocity is greater than a certain value, namely $v_{cut,out}$ speed or lower than a certain value, namely $v_{cut,in}$ speed, then the wind turbine must be curtailed; it does not operate. Eqn. 2.24 formulates the operating power coefficient of a wind turbine when the velocity value is in between $v_{cut,in}$ and $v_{cut,out}$ (Dolatabadi et al., 2017). Please note that the $v_{rated}$ value is the wind speed when the turbines operate at maximum power and the operating power coefficient is multiplied by the rated power of the wind turbine to calculate the generation as given by Eqn. 2.16.

Finally, the proposed MILP model accounts for the equipment efficiencies as constant values. For practical reasons, we did not include the effect of cycle efficiencies into the MILP model, whereas the structure of the model is compatible for including those cycle efficiencies for more detailed applications. Extra binary variables should be included to formulate the operation windows in that case. It should also be noted that the proposed model does not consider transmission costs, switch on-off costs for generators, ramping limits, and AC-DC (alternating current-direct current) conversion for the electricity. On the other hand, all these features are easily introducible to the decision-making model based on specific needs due to the generic structure.

## 3. Scenario Generation

The scenario generation approach determines the risk-awareness during the solution of the optimization problem. The variety of scenarios are expected to provide a more flexible architecture to address the needs at unlikely/unexpected cases. Our scenario generation method extends the scenario spectrum by propagating parameter uncertainty based on the statistical measures from the empirical model using the actual historical data.

A flowchart of the proposed scenario generation framework is shown in Fig 3.1. The historical data of the air temperature, solar irradiance, wind speed, population, and electricity consumption are obtained via an API (application programming interface) and distributors for the validation of the proposed empirical models (Dark Sky API, 2020). Next, the prediction covariance matrix and intervals are calculated from the propagation of model parameter uncertainty through Taylor series expansion, as in Equation 3.10. Training and parameter uncertainty propagation are done by MATLAB nonlinear regression toolbox using prediction confidence intervals (nlpredci) function, which calculates Jacobian and covariance matrix of parameters of the models with mean squared error (MSE) (Mathworks, 2022). Finally, scenarios are generated by inverse transform sampling at desired likelihood. The aforementioned steps will be explained in detail in the following subsections.

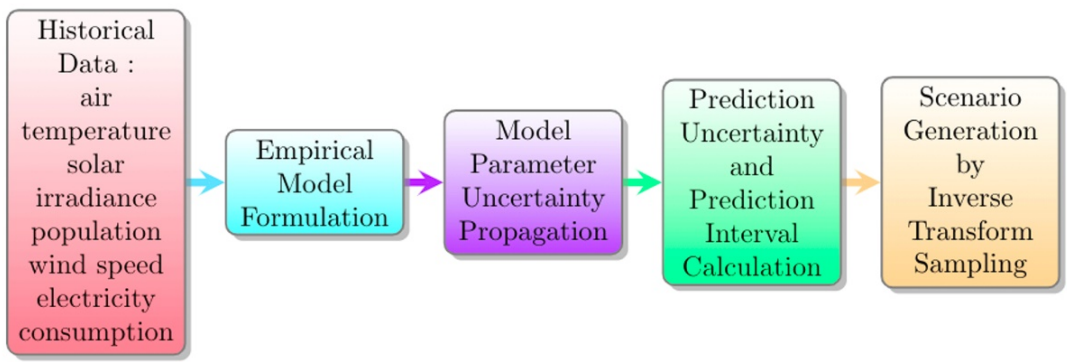

**Figure 3.1.** Scenario generation method flowchart.

Scenario reduction is a procedure to decrease the size of a scenario tree. If the branches of a scenario tree are not constrained or not predefined, they will aggregate, which in turn will increase computational complexity. In this case, scenario reduction is required to solve the corresponding stochastic problems in scalable times (Kaut and Wallace, 2003). Proposed method overcomes branch aggregation by providing a predetermined prediction interval. In other words, this method is able to decrease computational complexity since it does not need reduction after scenario tree generation.

## 3.1. Formulation of the Empirical Models

All empirical equations showing long or short-term behavior have flexible architectures to represent different profiles under distinct parameter values. As a result, the uncertainty region obtained from these models covers a significant fraction of historical data and would provide satisfactory prediction performance once the assumptions are still valid in the future. Therefore, the equations may be expected to address the major considerations such as trend, volatility, and conditioning (Rios et al., 2015; Sildir et al., 2021). Corresponding empirical model parameters, with index *l*, are obtained from:

$$\min_{p_l} \sum_{m=1}^{N} \| f_l(t_m, u_m, p_l) - y_m \| \tag{3.1}$$

where $N$ is the number of measurement samples; $p_l$ is the vector of parameters in equation $l$; $t_m$ is the time instant of measurement; $u_m$ is the input to the equation; $y_m$ is the measurement; $f_l$ is the proposed empirical equation for the prediction of $y_m$. Time variable "$t_m$" is continuous for all empirical formulations formed in this study.

### 3.1.1. Solar Irradiance

Sahin et *al*. developed the monthly solar irradiance model for 73 locations, including Yalova, through both Artificial Neural Networks (ANNs) and Multiple Linear Regression (MLR) using three years of training data (Şahin et al., 2013). Ozgoren et *al*. used a higher number of inputs, including air temperature, in order to train ANNs and develop nonlinear empirical equations (Ozgoren et al., 2012). In a different way, Bulut and Buyukalaca generated a trigonometric method which considers the day of the year as additional input and predicted the global-solar radiation with high accuracy based on 15 years of training data (Bulut and Büyükalaca, 2007). Sozen et *al*. trained ANNs based on the location, month, mean temperature, and sunshine duration (Sözen et al., 2005). There are also many empirical solar radiation models available (Arslanoglu, 2016; Ayvazoğluyüksel and Filik, 2017).

In this study, similar to the aforementioned studies above, historical data showed that solar irradiance at a particular time is highly dependent on time and temperature, requiring additional effort for air temperature forecast. The mean daily temperature is calculated from:

$$T_{daily} = p_{daily,1} \cdot \sin\left(\frac{2\pi}{365}(t + p_{daily,2})\right) + p_{daily,3} \tag{3.2}$$

where $T_{daily}$ is the daily mean air temperature; $t$ is the time; $p_{daily,1}$, $p_{daily,2}$, and $p_{daily,3}$ are the empirical model parameters to be estimated from the historical data; $\pi$ is 180°. Once $T_{daily}$ is calculated from Eqn. 3.2, the momentary air temperature value is obtained from:

$$T_{air} = \left(p_{airT,1} \cdot t_h^4 + p_{airT,2} \cdot t_h^3 + p_{airT,3} \cdot t_h^2 + p_{airT,4} \cdot t_h + p_{airT,5}\right) \cdot T_{daily} \tag{3.3}$$

where $p_{airT,1}$, $p_{airT,2}$, $p_{airT,3}$, $p_{airT,4}$, and $p_{airT,5}$ are the empirical model parameters; $t_h$ is the hour of the day. Eqn. 3.3 maps the daily mean temperature to minute-based predictions with a quite flexible architecture to represent various daily temperature profiles once the parameters have a certain variation range.

Once the air temperature at a particular time is calculated, the solar irradiance ($S$) is given by:

$$S = \frac{T_{daily}}{\sigma_s . \sqrt{(2\pi)}} e^{-\left(\frac{1}{2}\right)\left(\frac{t_h - \mu_s}{\sigma_s}\right)^2} \tag{3.4}$$

where $\sigma_s$ (standard deviation of solar irradiance) and $\mu_s$ (mean of solar irradiance) are:

$$\mu_s = p_{solar,\mu 1} . sin\left(\frac{2\pi}{365}(t + p_{solar,\mu 2})\right) + p_{solar,\mu 3} \tag{3.5}$$

$$\sigma_s = p_{solar,\sigma 1} . sin\left(\frac{2\pi}{365}(t + p_{solar,\sigma 2})\right) + p_{solar,\sigma 3} \tag{3.6}$$

$p_{solar,\mu 1}$, $p_{solar,\mu 2}$, $p_{solar,\sigma 1}$, and $p_{solar,\sigma 2}$ are the empirical model parameters. The variables in Eqns. 3.5 and 3.6 modify the Gaussian type behavior in Eqn. 3.4 by changing the sun peak time and the duration of sunlight availability, which is longer in hot seasons.

### 3.1.2. Wind Speed

Weibull and Rayleigh distribution are two important approaches to estimate the wind speed distribution (Gungor et al., 2020). Weibull parameters are reported in the literature for different locations (Akgül and Şenoğlu, 2019; Oner et al., 2013; Usta, 2016). Two-parameter Weibull distribution, Rayleigh distribution, and Inverse Weibull distribution for six cities, including Yalova, were studied (Dokur et al., 2019). In this study, wind speed is estimated by the traditional Weibull distribution type equation.

### 3.1.3. Electricity Demand

Electricity demand is a function of population and energy consumption per capita, in addition to seasonal variations and other uncertainties such as economic conditions, electricity price, potential for growth of industrial sectors (Hamzaçebi et al., 2019). In this study, the total electricity demand for a particular ($E_{daily}$) is calculated from:

$$E_{daily} = p_{e,daily,1} . sin\left(\frac{2\pi}{365}(t + p_{e,daily,2})\right) + p_{e,daily,3} + p_{e,daily,4} . P_t \tag{3.7}$$

where $P_t$ is the population at time $t$; $p_{e,daily,1}, p_{e,daily,2}, p_{e,daily,3}$, and $p_{e,daily,4}$ are the empirical model parameters. Eqn. 3.7 considers both the seasonal variation due to sine function and timely increase due to $P_t$. Note that the electricity consumption per capita is assumed constant. $P_t$ is calculated from:

$$P_t = P_0 + P_{p,1} \cdot t + P_{p,2} \cdot t^2 \tag{3.8}$$

where $P_0$ is the reference population; $P_{p,1}$ and $P_{p,2}$ are the empirical model parameters. Having calculated the daily electricity demand from Eqn. 3.7, the daily distribution can be obtained from:

$$E_h = \left(p_{E,h,1} t_h^4 + p_{E,h,2} t_h^3 + p_{E,h,3} t_h^2 + p_{E,h,4} t_h\right) \cdot E_{daily} \tag{3.9}$$

where $E_h$ is the energy demand at a particular time; $p_{E,h,1}$, $p_{E,h,2}$, $p_{E,h,3}$, and $p_{E,h,4}$ are the empirical model parameters for hourly energy demand. Eqn. 3.7 and Eqn. 3.9 show similarity with the results in the literature that hourly or seasonal electricity consumption in Turkey exhibited a sinusoidal trend (Bassols et al., 2002; Li et al., 2006; Preuster et al., 2017).

## 3.2. Uncertainty Propagation

The solution of the nonlinear optimization problem in Eqn. 3.1 provides a set of mean parameters. The resulting parameter values depend on the mathematical architecture of the empirical model in addition to data, noise, and optimization algorithm. Once the number of parameters is high for a particular model, some parameters might have a relatively lower impact on the outputs. Furthermore, different parameters might have a similar impact on the output due to strong correlations among them, causing identifiability problems (Günay, 2016; Yukseltan et al., 2017).

Parameter identifiability is evaluated through parameter covariance matrix (Schittkowski, 2007; Yetis and Jamshidi, 2014). This matrix provides a mathematical insight into variations from the mean parameter values. Different parameter samples can be obtained using this covariance matrix through Cholesky decomposition-based methods (Mclean and Mcauley, 2012), and such samples are likely to result in a different prediction profile. Thus, the probable variations in the parameters introduce the alternative profiles to represent actual uncertainty.

The difference in the prediction profile is quantified by both the change in the parameter and the output sensitivity to parameter. There are many methods to introduce the impact of the parameter uncertainty on the model predictions (Briggs et al., 2012; Kuchuk et al., 2010). One of the methods (Chakraborty, 2006) utilizes the Taylor series approximation to propagate the parameter uncertainty to the outputs using:

$$cov_y = J cov_p J^T \tag{3.10}$$

where $cov_y$ is the output covariance matrix; $cov_p$ is the parameter covariance matrix; $J$ is the Jacobian matrix evaluated at the desired prediction regime. Eqn. 3.10 takes only the first term of the Taylor, and higher accuracy is needed (Zouaoui and Wilson, 2001). In this particular case, $cov_y$ is a scalar, and in practice, represents the variance of the output variable.

The variance from Eqn. 3.10 defines the distribution of the mean prediction as a result of parameter uncertainty and does not take epistemic uncertainty sources (Kuczera and Parent, 1998), such as algorithmic uncertainty and structural uncertainty, into account. In order to consider the excluded sources, in addition to aleatoric uncertainties (Tellinghuisen, 2001), a lumped uncertainty term is additionally added to Eqn. 3.11, assuming normally distributed model prediction residuals. Such prediction uncertainty measure can easily be converted to the prediction intervals (Yegnan et al., 2002), which define the probability of a new measurement to exist between particular upper and lower bounds, under desired confidence level as shown in Eqn. 3.11 (Beven, 2016):

$$f_l(t_m, u_m, p_l) - \lambda\sqrt{cov_y + MSE} \leq y_m \leq f_l(t_m, u_m, p_l) + \lambda\sqrt{cov_y + MSE} \qquad (3.11)$$

where $\lambda$ is the coefficient of standard deviation of the mean value of $m^{th}$ output variable. It is related to the confidence level and determined by the inverse of Student's $t$ cumulative distribution function. *MSE* is the mean squared error.

After a probability distribution (or density) function of the output variables is obtained with the prediction interval calculated from Eqn. 3.11, inverse transform sampling method is used to generate the scenarios under their desired likelihood (Sugiyama, 2016). The inverse transform method enables sampling from a predefined probability range. If scenarios consist of only the most probable conditions, the decision-making model will tend to result in a very optimistic and fragile design. On the other hand, if scenarios include only the most improbable conditions, the corresponding design is conservative and suboptimal. Hence, in this study, the likely scenario is sampled around the mean prediction by one standard deviation interval, representing approximately 68.2 % of the output variable set. Also, the mid-likely and unlikely scenarios are sampled from standard deviations multiplied by two and three, respectively. Two standard deviation interval accounts for almost 27.2% of the output variable set while three standard deviation interval expresses nearly 4.2% of it.

Yalova is a relatively small but densely populated city in Turkey and has significant development potential with industrial and geopolitical advantages. A simplified map of the region is shown in Fig. 3.2, which includes the industrial zones in rectangles. Actual electricity consumption data are obtained from UEDAŞ (Uludağ Electricity Distribution Company), and corresponding scenario generation results are shown in Figure 3.3.

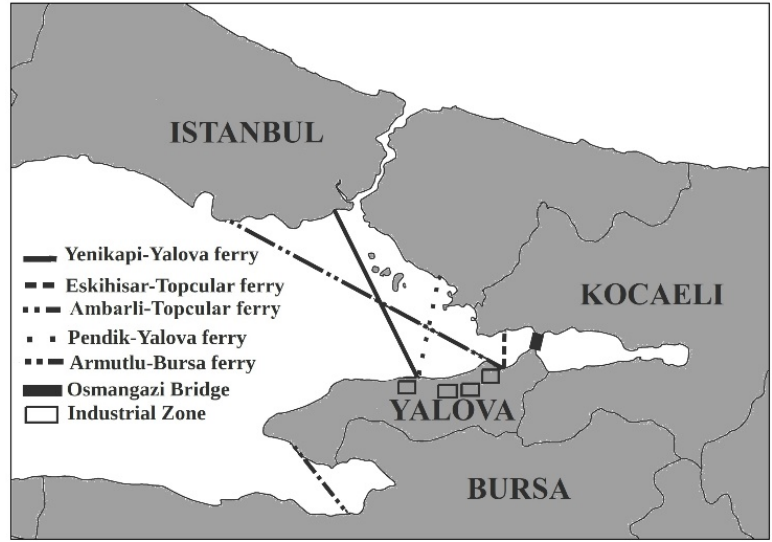

**Figure 3.2.** Yalova and surrounding region map.

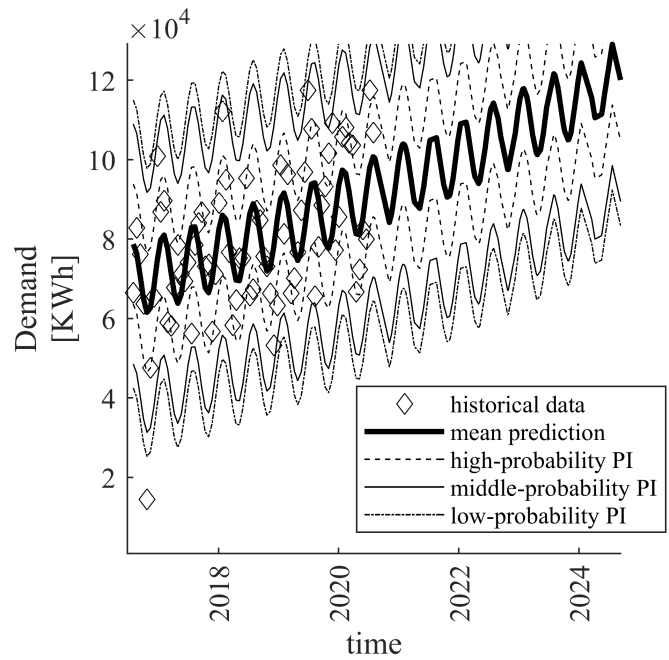

**Figure 3.3.** Energy consumption prediction intervals.

"Datenum" function in MATLAB is used to convert the datetime of the inputs to years. Hourly historical data are trained, and scenarios are generated over the next 20 years (from 2020 to 2040). In Figure 3.3, certain intervals are zoomed in to avoid a messy and complicated look. Only electricity consumption plot is drawn as an example demonstration of prediction intervals at different likely levels.

## 4. Results and Discussion

The proposed MILP model is formulated in GAMS (The General Algebraic System). CPLEX is the chosen MILP solver, which uses a branch and cut algorithm (CPLEX, 2022). In the solver options, ''optCA'' and ''optCR'', which are absolute and relative optimality criteria, respectively, are set to $10^{-2}$.

The deterministic MILP optimization model is formulated to provide one-third of the total electricity demand of Yalova. The remaining demand is assumed to be met by other private industry contracts and newly installed nuclear reactors among the nation. The lifespan of the project is twenty years. Ten different days, each of which is further divided into 30-minute intervals, are selected from the 20 years of forecast regime, and scenarios are generated for the corresponding instants. A day is also represented with forty-eight time intervals. Those time instants and scenarios are assumed to cover a wide range of requirements for the lifetime. However, a higher number of representative days could also be selected, which would, on the other hand, increases the computational complexity of the MILP model. Minimizing net present cost is determined as the objective function and the number of total equipment is limited to ten due to area restrictions.

Every piece of equipment has a different $CO_2$ footprint. The $CO_2$ emission ratios are collected from the same resources in the literature that the equipment's fuel consumption and energy generation ratios are obtained (Appendix Table A.3 and Table A.4). The $CO_2$ limit in the generation table in the appendix represents the emission rates produced from a piece of equipment when it runs at the unit power (kW) for unit time. A wide range of candidate equipment is introduced to the decision-making model, and the superstructure of the proposed energy system is shown in Figure 4.3.a. This candidate equipment consists of three main types, namely batteries, generators, and renewable equipment. The renewable energy equipment is wind turbines, solar cells, and biomass generators which use natural resources like wind, solar irradiance, and biomass, respectively. In addition, bio-gasifier can also be a renewable equipment using woody biomass as raw material. Electrolyzers and methanation reactors are also added as candidate equipment in a compact form. They are coupled in a way that the output hydrogen that comes from the electrolyzer is the input to the methanation reactor to produce synthetic natural gas. Thus, hydrogen becomes an intermediate product removing any storage related problems. Electricity is also needed for hydrogen production through electrolyzers. The data including consumption, generation, cost and operation limits of resources and candidate equipment, are given in the Appendix (Tables A.1- A.4.). Optimal results for three different scenarios are obtained using actual data and the proposed scenario generation method.

We generated three different scenarios, namely likely, mid-likely and unlikely. In each scenario, 39 candidate equipment are introduced to the decision-making model. The chosen equipment for different scenarios are given in Figure 4.3.b. One of the candidate equipment sets is the IGCC's, namely integrated gasification combined cycle (Wang, 2017). While the IGCC-2 has the same working principle with the IGCC-1 with doubled capacity, the IGCC-3 includes a water gas shift reactor and $CO_2$ capture unit. This is why the IGCC-3 has a relatively low $CO_2$ emission rate compared to the IGCC-1 and IGCC-2 (*Updated Capital Cost Estimates for Electricity Generation Plants*, 2010). However, IGCC-3 is not chosen.

Optimal solutions show that, for all cases, the decision-making model prefers to buy electricity from the maxi-grid up to its maximum level (3MWh) except for a few time intervals since purchasing electricity does not require extra expenses such as equipment installation cost or maintenance costs. On the other hand, it is not preferable to buy electricity when there is surplus electricity inside the grid. Since purchasable electricity from the national grid is limited due to stability issues, the rest of the demand should be met inside the grid.

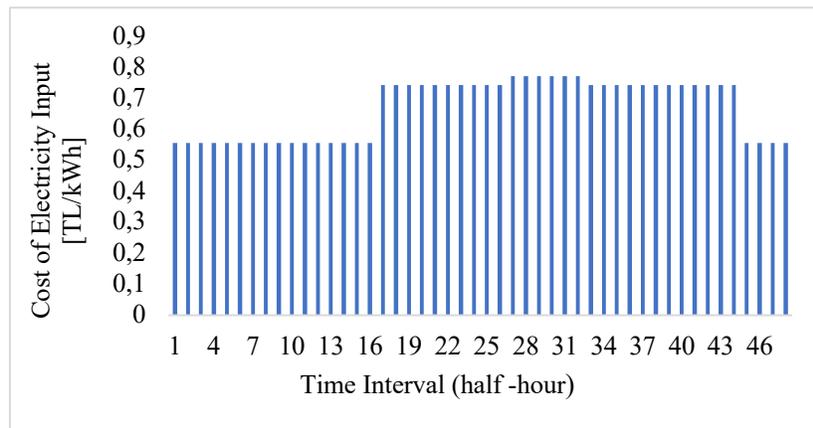

**Figure 4.1.** Cost for unit system input of electricity (taken from EPDK).

Electricity demand profiles within a day for three different scenarios are shown in Figure 4.2. Please note that Yalova is a heavily industrial city, which in turn contributes to higher demand profiles at night as well. Once the daily demand summations are compared, the likely scenario has the highest demand value while the mid-likely scenario has the lowest demand value.

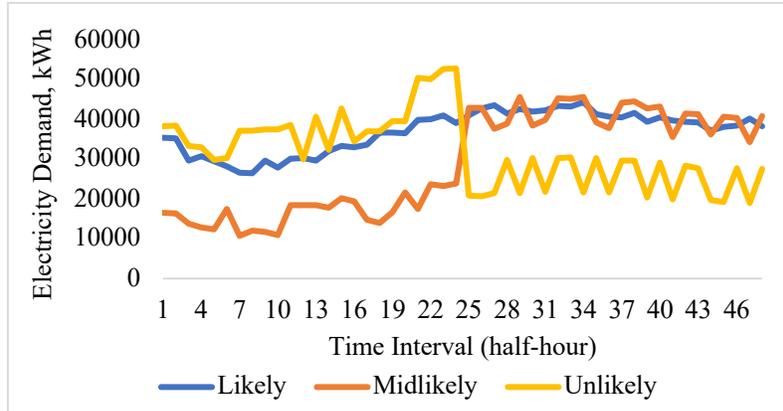

**Figure 4.2.** Electricity demand data for three scenarios**.**

The total number of chosen equipment is nine for the likely and mid-likely scenarios and ten for the unlikely scenario. The selection of the optimal equipment and corresponding sizes and rated powers are provided in Table 4.1. In optimal results, the electricity demand is mostly supplied with PV units for the likely and mid-likely scenarios and bio-gasifiers for the unlikely scenario. There is a contribution to power supply from wind turbines only for mid-likely scenario stemming from higher wind speed rates than other scenarios. Finally, the electricity input is 3MWh for most of the time intervals for three scenarios where it is the upper limit of electricity input from the national grid.

### 4.1. Power Balance
#### *4.1.1. Likely Scenario*
The optimal cost is found as 9.6 billion TRY when a 12% discount rate is considered. Only one PV unit is chosen, and this unit has the highest rated power; wind turbine is not chosen. Naturally, the operating power of the PV unit is directly related to the UV (Ultraviolet) light (or solar irradiance). Accordingly, Figure 4.4 shows that the electricity demand is mostly supplied with PV unit when there is available sunlight. The peaks in Figure 4.4, which belong to the PV equipment, overlap with the solar irradiance data (Figure 4.7).

**Table 4.1.** Optimal solutions of three scenarios in terms of equipment decision and sizes.

| Chosen Equipment (Likely Case) | Rated Power (kW) | Chosen Units (Mid-Likely Case) | Rated Power (kW) | Chosen Units (Un-Likely Case) | Rated Power (kW) |
|---|---|---|---|---|---|
| Photovoltaic-3 (PV) | 133308.2 | Wind Turbine-1 | 8904.2 | Photovoltaic-1 (PV) | 21901.6 |
| Oil Co-generator -2 | 7717.2 | Photovoltaic-2 (PV) | 97235.4 | Oil Co-generator -1 | 3675.6 |
| Oil Co-generator -3 | 3076.4 | Oil Co-generator -1 | 6527 | Oil Co-generator -2 | 7712.4 |
| Reciprocating ICE-2 | 1121 | Oil Co-generator -2 | 6459.5 | Oil Co-generator -3 | 7713 |
| Reciprocating ICE-3 | 9341 | Oil Co-generator -3 | 6314.9 | Reciprocating ICE-2 | 1121 |
| Micro Turbine-1 | 240 | Reciprocating ICE-2 | 1121 | Reciprocating ICE-3 | 9341 |
| Micro Turbine-2 | 320 | Reciprocating ICE-3 | 9341 | Micro Turbine-2 | 320 |
| Biogasifier-1 | 6600 | Micro Turbine-2 | 320 | Micro Turbine-3 | 950 |
| Biogasifier-2 | 11600 | Biogasifier-2 | 11600 | Biogasifier-1 | 6600 |
| | | | | Biogasifier-2 | 11600 |

Figure 4.4. also shows that there is excess production in the time intervals 26 and 28, each time interval corresponding to half an hour. Those intervals match exactly with the two intervals of solar irradiance data, the third-highest and the highest peak. The grid does not prefer to buy electricity from the national grid at these time intervals because there is surplus electricity. The amount of surplus electricity is 1266 and 2116 kWh for these time intervals. It should be mentioned here that the grid is let have little surplus electricity, which must be less than 5 % of the highest demand. This surplus can easily be managed via internal usage or transfer to the maxi-grid. Eventually, the MILP model does not decide to install any storage device due to the installation, maintenance costs and demand level.

As aforementioned before, one of the popular methods to store excess electricity produced from renewable energy sources is to convert the surplus electricity into chemical products, e.g., methane or ammonia. The idea is to convert that synthesized natural gas into electricity back again when electricity prices or importing natural gas are high. This compact formulation would be a feasible choice instead of storing the surplus electricity in a battery. However, electrolyzers and methanation reactors are not chosen in the optimal solutions of the scenarios. We observed that investing in this compact form is simply not economically favorable due to the installation costs and limited renewable resources such as wind speed and solar irradiance. The model decides to import as much electricity outside the grid as possible and meet the electricity demand without converting excess electricity to chemical intermediate products.

Sensitivity analysis could also be performed to detect the conditions where green hydrogen and synthetic methanation become economically feasible. To verify the selections, we reduced the costs of installing the methanation reactor, electrolyzer and we allowed purchasing hydrogen and water outside the grid with reduced costs. This way, the model indeed chooses to install the methanation reactor. As a result, the opti-mal solution suggests importing hydrogen and water from outside the grid to operate the methanation reactor. At the same time, the reciprocating engine uses the synthesized natural gas when there is peak electricity demand, reducing overall electricity purchasing prices significantly.

On the other hand, when purchasing the hydrogen from outside is not allowed, or when the hydrogen import price is set to be the real value, installing the compact electrolyzer-methanation pair is not preferred. Consequently, this concept is not optimal for Yalova due to the installation and raw material prices and current demand levels. Nevertheless, this result can only be confirmed under the current price and demand conditions for Yalova. Different outcomes might be observed for larger cities and regions, reduced costs, changing demand profiles and tightening emission limits and regulations.

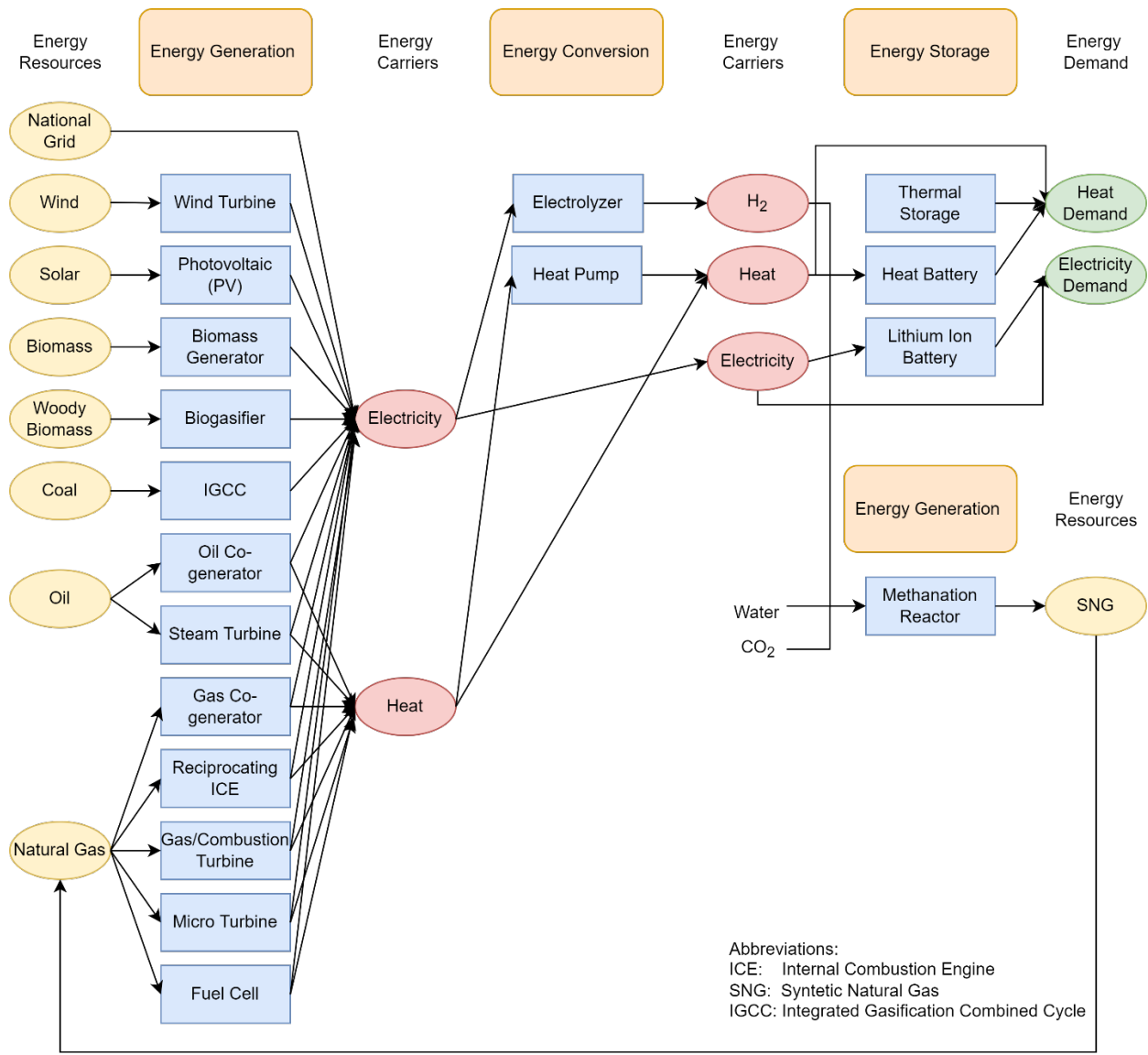

**Figure 4.3.a.** Superstructure of the energy system. (Zhou et al., 2013)

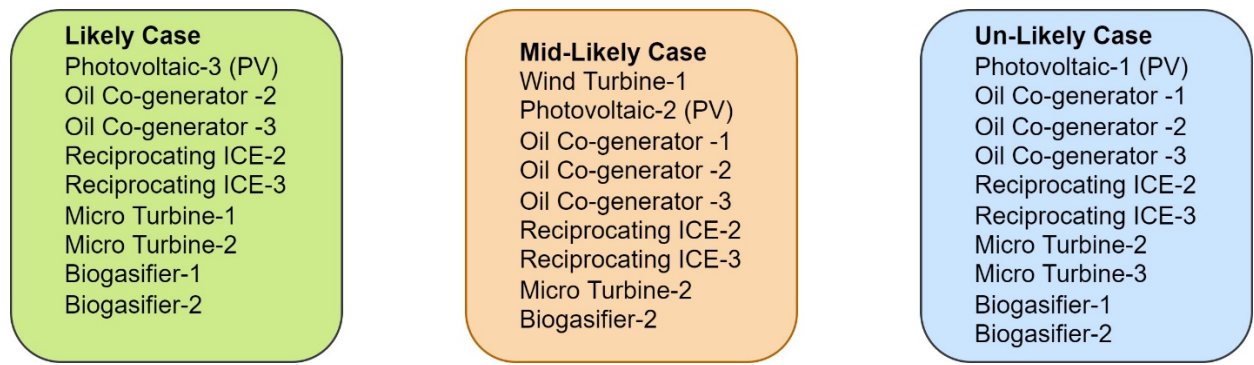

**Figure 4.3.b.** Chosen equipment for different cases among candidate equipment.

*4.1.2. Mid-Likely Scenario*

The optimal equipment list is seen in Table 4.1. for the mid-likely scenario, with the corresponding rated power values. The optimal cost is found as 8.01 billion TRY.

It is seen from Figure 4.6. that the PV unit has the highest rated power. Unlike the likely scenario, there is also a contribution to electricity supply from wind turbines as renewable equipment in this scenario. This is due to the fact that the mid-likely scenario includes higher wind speed rates. Mean and maximum wind speed rates for three different scenarios are provided in Table 4.2. and Table 4.3. Wind speed data is given in Figure 4.5. There is an operational constraint that wind speed should be greater than 2.5 [m/sec] in order to operate. Hence, the values lower than the 2.5 [m/sec] show insufficient wind speed values. The time intervals that wind turbine operates in the mid-likely scenario match with the wind speed peaks in Figure 4.5.

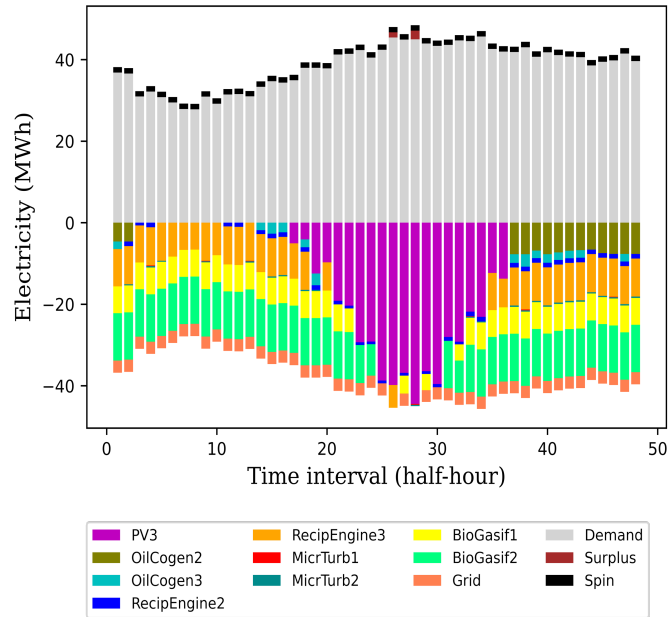

**Figure 4.4.** Optimal deterministic results for a day in the 1$^{st}$ year for likely scenario.

**Table 4.2.** Mean and maximum wind speed rates among all wind speed values [m/sec].

| Scenarios | Likely | Mid-Likely | Unlikely |
|---|---|---|---|
| Mean Value | 2.52 | 2.39 | 2.15 |
| Maximum Value | 5.14 | 5.87 | 4.51 |

**Table 4.3.** Mean and maximum wind speed rates among wind speed values greater than 2.5 [m/sec].

| Scenarios | Likely | Mid-Likely | Unlikely |
|---|---|---|---|
| Mean Value | 3.44 | 3.72 | 3.34 |
| Maximum Value | 5.14 | 5.87 | 4.51 |

However, the contribution of the wind turbine is relatively small compared to the solar cell. Electricity demand is mostly supplied with the bio-gasifier and reciprocating engine units after the solar cell, followed by oil co-generators. The initial and maintenance costs of the oil co-generators are lower than the reciprocating engines and bio-gasifier. One of the reasons that the bio-gasifier and reciprocating engines are operated more frequently than the oil co-generators is because they process cheaper raw materials than oil. Another reason is that the reciprocating engines and the bio-gasifier have lower $CO_2$ emission rates than the oil co-generators.

Figure 4.6. shows the scheduling decisions for the mid-likely scenario for the first year. There is a small amount of surplus electricity produced in the 21$^{st}$ time interval. In this time interval, there is a high peak in the solar irradiance data (Figure 4.7). In this interval, the model does not buy electricity because there is already surplus electricity in the grid, albeit very small.

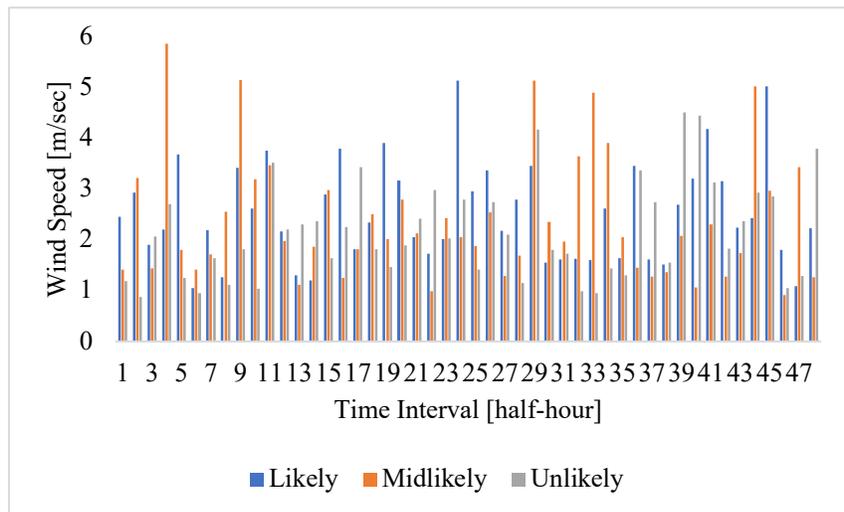

**Figure 4.5.** Wind Speed data for three scenarios.

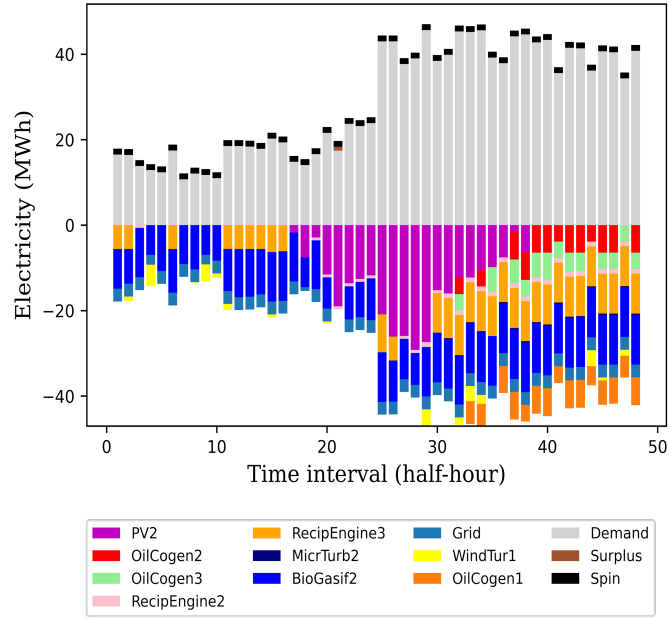

**Figure 4.6.** Optimal deterministic results for a day in the 1$^{st}$ year for mid-likely scenario.

Similar to the likely scenario, obtaining surplus electricity shows that the renewable equipment is not sufficient to operate the compact electrolyzer-methanation pair even if the initial cost of this compact pair is set to zero.

### 4.1.3. Un-likely Scenario

The chosen equipment and corresponding rated powers for the unlikely scenario are given in Table 4.1. The optimal cost is found as 7.36 billion TRY. The electricity demand is mostly supplied with the bio-gasifier unit in the unlikely scenario case, even though the rated power of solar cells is greater than the sum of the rated powers of bio-gasifiers. However, bio gasifiers operate in every time interval, mostly at their maximum rated power. Hence, the total operating power of the bio-gasifiers is greater than the solar cell for the day in year one, as shown in Figure 4.8. It is followed by reciprocating engines and oil co-generators, respectively.

Solar irradiance data graph for every scenario is given in Figure 4. 7. Typically, the operating power of the units is directly proportional to their rated powers, which are reported in Table 4.1. Please note that the rated power of the photovoltaic generator in the unlikely scenario is the lowest among all scenarios. Therefore, while power supply contribution from the PV unit is the highest in the likely scenario, it is the lowest in the unlikely scenario.

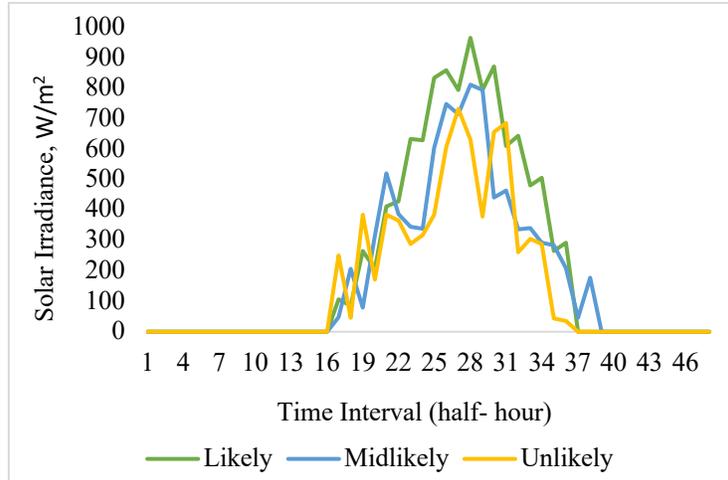

**Figure 4.7** Solar irradiance data graph for every scenario.

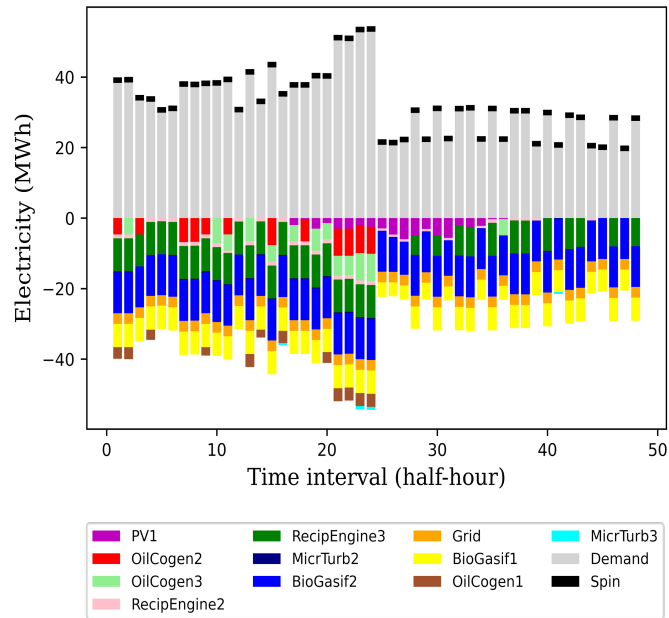

**Figure 4.8.** Optimal deterministic results for a day in the 1$^{st}$ year for unlikely scenario.

## 4.2. Heat Balance

Heat purchase from outside of the grid is not allowed. For brevity, we provide only the heat balance profile for the likely scenario in Figure 4.9. As shown in Figure 4.9, excess heat is generated by the generators inside the system. Those generators are combined heat and power units, and they both produce heat and electricity at certain ratios, which are usually greater than one, as given in Appendix Table A.4. The main aim of the proposed model is to meet the electricity demand. Therefore, it is inevitable to produce surplus heat while meeting the electricity demand. One way to tackle with surplus heat is to use equipment that converts thermal energy to electricity. Yet, these units are usually expensive and exhibit very low cycle efficiencies and therefore are not included in the candidate pool in this work.

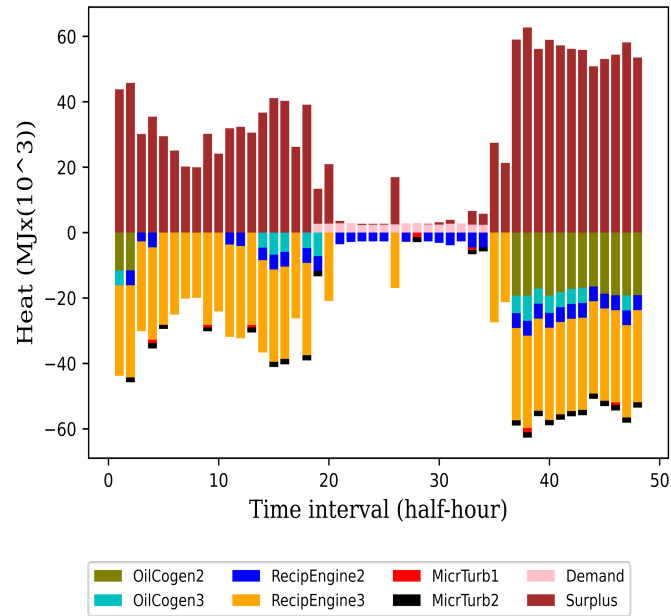

**Figure 4.9.** Heat balance for a day in the 1st year for likely scenario.

**Table 4.4.** Sizes and CPU times.

|  | Likely | Mid-Likely | Unlikely |
|---|---|---|---|
| CPU Time (sec) | 1763,82 | 3613.9 | 11469.68 |
| Continuous Variables | 337,169 | 337,169 | 337,169 |
| Discrete Variables | 31,719 | 31,719 | 31,719 |
| Number of Constraints | 408,311 | 408,311 | 408,311 |

In general, it is fair to state that similar equipment sets are proposed for three different scenarios. Bio-gasifier equipment is chosen for all scenarios. The main reason is the usage of cheaper raw material, woody biomass. Turkey, especially Marmara and Aegean regions, are known to have suitable conditions for raising such plants quickly. Similarly, current trends in Turkey show that biomass is and will be a nice candidate for the transition to the zero-carbon power industry (Difiglio et al., 2020). Finally, please note that the required CPU (central processing unit) times needed to solve for three scenarios are significantly different even though the model sizes are exactly the same, as shown in Table 4.4.

Other equipment commonly chosen for all cases are reciprocating internal combustion engines, oil co-generators and micro turbines. Reciprocating ICE's and micro turbines use natural gas as raw material, whereas oil co-generator uses oil. Natural gas is slightly expensive than oil. On the other hand, the energy content of natural gas is greater than oil, so that the electricity production rate of natural gas is also greater than oil ( Oil and Gas Industry Overview, 2021).

Finally, in order to satisfy the emission constraints, renewable equipment is more favorable for the decision-making model due to their low emission values. To observe that, sensitivity analysis is made for each scenario to see which equipment is chosen when the emission limit is relaxed or tightened. The optimal solutions tend to prefer the wind and solar units with higher rated power values when the emission limit is tightened even though they have higher initial costs than the non-renewable equipment. When the emission limit is relaxed, the decision-making model prefers primarily non-renewable equipment due to their lower initial costs. However, relaxed $CO_2$ emission limits are not the expected future trend nowadays. On the contrary, $CO_2$ cap trade incomes are expected to justify the installation of renewable equipment and should also be considered for future studies. We also believe that the aforementioned analysis will justify the requirement of green hydrogen for net-zero emission targets for different locations and changing trends due to the Paris agreement. Furthermore, extended formulations such as stochastic programming combined with the suggested holistic scenario generation and optimization model might provide more reliable decisions as climate shift regulations, such as net-zero emission targets and carbon dioxide taxing and cap trade, have been changing rapidly.

## 5. Conclusion

In this study, a generic MILP decision-making model is developed for the optimal integration of renewable and non-renewable units while considering weather and demand uncertainties and chemical conversion such as green hydrogen and synthetic methanation. The case studies shown in this paper are related to the industrial city Yalova in Turkey. A wide range of equipment is selected as the candidate equipment, including non-renewable and renewable generators, electrolyzers, compact Power-to-X concept plants, green hydrogen, and storage devices. In addition to the proposed MILP decision-making model, this work also suggests a practical scenario generation method to produce dependable scenarios using actual data. Accordingly, three different scenarios, namely likely, mid-likely and unlikely, are generated. Afterward, these generated scenarios are analyzed separately.

Optimal results for different scenarios show that the electricity demand is mostly supplied with PV units for the likely and mid-likely scenarios and bio-gasifiers for the unlikely scenario. Moreover, the total cost, including the initial investment, operation and maintenance costs, is 9.6, 8.01, and 7.36 billion TRY for the likely, mid-likely, and unlikely scenarios, respectively, when a 12% inflation rate is considered. It is also observed that the $CO_2$ emission constraint has a vital impact on optimal decision-making. In that sense, bio-gasifier and PV units are the ones contributing the most to the production of low-carbon power. It can be stated that the optimal solution of the model will include mainly the renewable energy-based equipment when the $CO_2$ emission rates are tightened, which is obtained through sensitivity analysis. On the other hand, green hydrogen and synthetic methanation concepts are not economically preferable for the investigated location currently, when optimal design and operation criteria are assessed simultaneously, unlike many other contributions in the literature. Nevertheless, since the emission limits are getting tightened year by year with international regulations, it is expected that constructing renewable resource-based grids would be more preferable for Turkey in the future, as sanctions will be tighter. Thus, earlier investment decisions must be discussed for the investigated location. Furthermore, including the cap and trade of $CO_2$ as income streams can also reduce the overall costs significantly. Finally, the generic nature of the proposed model makes it suitable for two-stage stochastic formulations, which is the current research focus of the authors.

**Acknowledgments**

This publication has been produced benefiting from the 2232 International Fellowship for Outstanding Researchers Program of TUBITAK (Project No: 118C245). However, the entire responsibility of the publication belongs to the owner of the publication. We thank Uludag Electricity Distribution Company (UEDAS) for sharing the actual electricity consumption data for Yalova.

# APPENDIX

## A. Tables

**Table A.1.** Power and capacity limits (Darrow et al., n.d.; EIA, n.d.; EPA and CHP, n.d.; *Updated Capital Cost Estimates for Electricity Generation Plants*, 2010).

|  | Minimum Rated Power [kW] | Maximum Rated Power [kW] | Minimum Capacity | Maximum Capacity |
|---|---|---|---|---|
| Wind Turbine-1 | 100 | 200000 | 0 | 0 |
| Wind Turbine-2 | 100 | 200000 | 0 | 0 |
| Wind Turbine-3 | 100 | 200000 | 0 | 0 |
| Photovoltaic-1(PV) | 100 | 200000 | 0 | 0 |
| Photovoltaic-2(PV) | 100 | 200000 | 0 | 0 |
| Photovoltaic-3(PV) | 100 | 200000 | 0 | 0 |
| Biomass Generator | 2500 | 500000 | 0 | 0 |
| Gas Co-generator-1 | 3350 | 670000 | 0 | 0 |
| Gas Co-generator-2 | 3350 | 670000 | 0 | 0 |
| Gas Co-generator-3 | 3350 | 670000 | 0 | 0 |
| Oil Co-generator -1 | 750 | 150000 | 0 | 0 |
| Oil Co-generator -2 | 750 | 150000 | 0 | 0 |
| Oil Co-generator -3 | 750 | 150000 | 0 | 0 |
| Lithium-Ion Battery | 1000 | 500000 | 100000 kWh | 500000 kWh |
| Lithium-Ion Battery-1 | 1000 | 500000 | 100000 kWh | 500000 kWh |
| Heat Pump | 1000 | 300000 | 0 | 0 |
| Thermal Energy Storage | 1000 | 500000 | 30000 MJ | 120000 MJ |
| Reciprocating Internal Combustion Engine-1 | 100 | 100 | 0 | 0 |
| Reciprocating Internal Combustion Engine-2 | 1121 | 1121 | 0 | 0 |
| Reciprocating Internal Combustion Engine-3 | 9341 | 9341 | 0 | 0 |
| Gas/Combustion Turbine-1 | 3304 | 3304 | 0 | 0 |
| Gas/Combustion Turbine-2 | 7038 | 7038 | 0 | 0 |
| Gas/Combustion Turbine-3 | 9950 | 9950 | 0 | 0 |
| Steam Turbine-1 | 500 | 500 | 0 | 0 |
| Steam Turbine-2 | 3000 | 3000 | 0 | 0 |
| Steam Turbine-3 | 15000 | 15000 | 0 | 0 |
| Micro Turbine-1 | 240 | 240 | 0 | 0 |
| Micro Turbine-2 | 320 | 320 | 0 | 0 |
| Micro Turbine-3 | 950 | 950 | 0 | 0 |
| Fuel Cell-1 | 1400 | 1400 | 0 | 0 |
| Fuel Cell-2 | 400 | 400 | 0 | 0 |
| Fuel Cell-3 | 300 | 300 | 0 | 0 |
| Biogasifier-1 | 6600 | 6600 | 0 | 0 |
| Biogasifier-2 | 11600 | 11600 | 0 | 0 |
| IGCC-1 | 600000 | 600000 | 0 | 0 |
| IGCC-2 | 1200000 | 1200000 | 0 | 0 |
| IGCC-3 | 520000 | 520000 | 0 | 0 |
| Electrolyzer | 10000 | 100000 | 0 | 0 |
| Methanation Reactor | 51300 | 51300 | 0 | 0 |

Table A.2. Initial and maintenance costs for the first year (Darrow et al., n.d.; EIA, n.d.; EPA and CHP, n.d.; *Updated Capital Cost Estimates for Electricity Generation Plants*, 2010).

| Symbol | Initial Cost | | Maintenance Cost | |
|---|---|---|---|---|
| | Per Rated Power [TRY/kW] | Per Capacity [TRY/kWh-MJ] | Per Rated Power [TRY/kW] | Per Capacity [TRY/kWh-MJ] |
| Wind Turbine-1 | 20039.04 | 0 | 403.2 | 0 |
| Wind Turbine-2 | 20039.04 | 0 | 403.2 | 0 |
| Wind Turbine-3 | 20039.04 | 0 | 403.2 | 0 |
| Photovoltaic-1(PV) | 20744.64 | 0 | 248.6 | 0 |
| Photovoltaic-2(PV) | 20744.64 | 0 | 248.6 | 0 |
| Photovoltaic-3(PV) | 20744.64 | 0 | 248.6 | 0 |
| Biomass Generator | 28082.88 | 0 | 1814.4 | 0 |
| Gas Co-generator-1 | 8537.76 | 0 | 672 | 0 |
| Gas Co-generator-2 | 8537.76 | 0 | 672 | 0 |
| Gas Co-generator-3 | 8537.76 | 0 | 672 | 0 |
| Oil Co-generator -1 | 9172.8 | 0 | 530.9 | 0 |
| Oil Co-generator -2 | 9172.8 | 0 | 530.9 | 0 |
| Oil Co-generator -3 | 9172.8 | 0 | 530.9 | 0 |
| Lithium-Ion Battery | 0 | 12096 | 134.4 | 336 |
| Lithium-Ion Battery-1 | 0 | 12096 | 134.4 | 336 |
| Heat Pump | 16800 | 0 | 672 | 0 |
| Thermal Energy Storage | 0 | 268.8 | 0 | 26.9 |
| Reciprocating Internal Combustion Engine-1 | 20300 | 0 | 1512 | 0 |
| Reciprocating Internal Combustion Engine-2 | 16562 | 0 | 1149 | 0 |
| Reciprocating Internal Combustion Engine-3 | 10031 | 0 | 514.1 | 0 |
| Gas/Combustion Turbine-1 | 22967 | 0 | 762.05 | 0 |
| Gas/Combustion Turbine-2 | 14560 | 0 | 743.9 | 0 |
| Gas/Combustion Turbine-3 | 13832 | 0 | 725.6 | 0 |
| Steam Turbine-1 | 7952 | 0 | 604.8 | 0 |
| Steam Turbine-2 | 4774 | 0 | 544.32 | 0 |
| Steam Turbine-3 | 4662 | 0 | 362.9 | 0 |
| Micro Turbine-1 | 19040 | 0 | 729.1 | 0 |
| Micro Turbine-2 | 18060 | 0 | 592.2 | 0 |
| Micro Turbine-3 | 17500 | 0 | 725.8 | 0 |
| Fuel Cell-1 | 32200 | 0 | 2419.2 | 0 |
| Fuel Cell-2 | 49000 | 0 | 2177.3 | 0 |
| Fuel Cell-3 | 70000 | 0 | 2721.6 | 0 |
| Biogasifier-1 | 34392 | 0 | 2693.4 | 0 |
| Biogasifier-2 | 28646 | 0 | 2038.2 | 0 |
| IGCC-1 | 24955 | 0 | 830.11 | 0 |
| IGCC-2 | 22547 | 0 | 757.8 | 0 |
| IGCC-3 | 37450 | 0 | 971.4 | 0 |
| Electrolyzer | 5355 | 0 | 289 | 0 |
| Methanation Reactor | 5580.1 | 0 | 210.6 | 0 |

**Table A.3.** Consumption table at unit power [kW] in unit time [h] (Darrow et al., n.d.; EIA, n.d.; EPA and CHP, n.d.; *Updated Capital Cost Estimates for Electricity Generation Plants*, 2010).

| Equipment Name | Consumption | | | | | | | | | |
|---|---|---|---|---|---|---|---|---|---|---|
| | Elect. [kWh/kW/h] | Heat [MJ/kW/h] | Biomass [MJ/kW/h] | Gas [kWh/kW/h] | Oil [MJ/kW/h] | $CO_2$ [$gCO_2$/kW/h] | Wood F. [MJ/kW/h] | Coal [MJ/kW/h] | $H_2$ [kWh/kW/h] | Water [kg/kW/h] |
| Wind Turbine-1 | 0 | 0 | 0 | 0 | 0 | 0 | 0 | 0 | 0 | 0 |
| Wind Turbine-2 | 0 | 0 | 0 | 0 | 0 | 0 | 0 | 0 | 0 | 0 |
| Wind Turbine-3 | 0 | 0 | 0 | 0 | 0 | 0 | 0 | 0 | 0 | 0 |
| Photovoltaic-1(PV) | 0 | 0 | 0 | 0 | 0 | 0 | 0 | 0 | 0 | 0 |
| Photovoltaic-2(PV) | 0 | 0 | 0 | 0 | 0 | 0 | 0 | 0 | 0 | 0 |
| Photovoltaic-3(PV) | 0 | 0 | 0 | 0 | 0 | 0 | 0 | 0 | 0 | 0 |
| Biomass Generator | 0 | 0 | 20.07 | 0 | 0 | 0 | 0 | 0 | 0 | 0 |
| Gas Co-generator-1 | 0 | 0 | 0 | 3.81 | 0 | 0 | 0 | 0 | 0 | 0 |
| Gas Co-generator-2 | 0 | 0 | 0 | 3.81 | 0 | 0 | 0 | 0 | 0 | 0 |
| Gas Co-generator-3 | 0 | 0 | 0 | 3.81 | 0 | 0 | 0 | 0 | 0 | 0 |
| Oil Co-generator -1 | 0 | 0 | 0 | 0 | 11.3 | 0 | 0 | 0 | 0 | 0 |
| Oil Co-generator -2 | 0 | 0 | 0 | 0 | 11.3 | 0 | 0 | 0 | 0 | 0 |
| Oil Co-generator -3 | 0 | 0 | 0 | 0 | 11.3 | 0 | 0 | 0 | 0 | 0 |
| Lithium Ion Battery | 1.05 | 0 | 0 | 0 | 0 | 0 | 0 | 0 | 0 | 0 |
| Lithium Ion Battery-1 | 1.05 | 0 | 0 | 0 | 0 | 0 | 0 | 0 | 0 | 0 |
| Heat Pump | 1 | 0 | 0 | 0 | 0 | 0 | 0 | 0 | 0 | 0 |
| Thermal Energy Storage | 0 | 1.2 | 0 | 0 | 0 | 0 | 0 | 0 | 0 | 0 |
| Reciprocating Internal Combustion Engine-1 | 0 | 0 | 0 | 3.69 | 0 | 0 | 0 | 0 | 0 | 0 |
| Reciprocating Internal Combustion Engine-2 | 0 | 0 | 0 | 2.72 | 0 | 0 | 0 | 0 | 0 | 0 |
| Reciprocating Internal Combustion Engine-3 | 0 | 0 | 0 | 2.4 | 0 | 0 | 0 | 0 | 0 | 0 |
| Gas/Combustion Turbine-1 | 0 | 0 | 0 | 4.18 | 0 | 0 | 0 | 0 | 0 | 0 |
| Gas/Combustion Turbine-2 | 0 | 0 | 0 | 3.47 | 0 | 0 | 0 | 0 | 0 | 0 |
| Gas/Combustion Turbine-3 | 0 | 0 | 0 | 3.66 | 0 | 0 | 0 | 0 | 0 | 0 |
| Steam Turbine-1 | 0 | 0 | 0 | 0 | 57.4 | 0 | 0 | 0 | 0 | 0 |
| Steam Turbine-2 | 0 | 0 | 0 | 0 | 73.2 | 0 | 0 | 0 | 0 | 0 |
| Steam Turbine-3 | 0 | 0 | 0 | 0 | 49.2 | 0 | 0 | 0 | 0 | 0 |
| Micro Turbine-1 | 0 | 0 | 0 | 3.85 | 0 | 0 | 0 | 0 | 0 | 0 |
| Micro Turbine-2 | 0 | 0 | 0 | 3.57 | 0 | 0 | 0 | 0 | 0 | 0 |
| Micro Turbine-3 | 0 | 0 | 0 | 3.76 | 0 | 0 | 0 | 0 | 0 | 0 |
| Fuel Cell-1 | 0 | 0 | 0 | 2.35 | 0 | 0 | 0 | 0 | 0 | 0 |
| Fuel Cell-2 | 0 | 0 | 0 | 2.92 | 0 | 0 | 0 | 0 | 0 | 0 |
| Fuel Cell-3 | 0 | 0 | 0 | 2.13 | 0 | 0 | 0 | 0 | 0 | 0 |
| Biogasifier-1 | 0 | 0 | 0 | 0 | 0 | 0 | 20.5 | 0 | 0 | 0 |
| Biogasifier-2 | 0 | 0 | 0 | 0 | 0 | 0 | 20.49 | 0 | 0 | 0 |
| IGCC-1 | 0 | 0 | 0 | 0 | 0 | 0 | 0 | 9.18 | 0 | 0 |
| IGCC-2 | 0 | 0 | 0 | 0 | 0 | 0 | 0 | 9.18 | 0 | 0 |
| IGCC-3 | 0 | 0 | 0 | 0 | 0 | 0 | 0 | 11.3 | 0 | 0 |
| Electrolyzer | 1.32 | 0 | 0 | 0 | 0 | 0 | 0 | 0 | 0 | 0 |
| Methanation Reactor | 0 | 0 | 0 | 0 | 0 | 177.14 | 0 | 0 | 1.28 | 0.748 |

**Table A.4.** Generation table at unit power [kW] in unit time [h] (Darrow et al., n.d.; EIA, n.d.; EPA and CHP, n.d.; *Updated Capital Cost Estimates for Electricity Generation Plants*, 2010).

| Equipment Name | Generation | | | | | | | | | |
|---|---|---|---|---|---|---|---|---|---|---|
| | Elect. [kWh/kW/h] | Heat [MJ/kW/h] | Biomass [MJ/kW/h] | Gas [kWh/kW/h] | Oil [MJ/kW/h] | $CO_2$ [g$CO_2$/kW/h] | Wood F. [MJ/kW/h] | Coal [MJ/kW/h] | $H_2$ [kWh/kW/h] | Water [kg/kW/h] |
| Wind Turbine-1 | 1 | 0 | 0 | 0 | 0 | 0 | 0 | 0 | 0 | 0 |
| Wind Turbine-2 | 1 | 0 | 0 | 0 | 0 | 0 | 0 | 0 | 0 | 0 |
| Wind Turbine-3 | 1 | 0 | 0 | 0 | 0 | 0 | 0 | 0 | 0 | 0 |
| Photovoltaic-1(PV) | 1 | 0 | 0 | 0 | 0 | 0 | 0 | 0 | 0 | 0 |
| Photovoltaic-2(PV) | 1 | 0 | 0 | 0 | 0 | 0 | 0 | 0 | 0 | 0 |
| Photovoltaic-3(PV) | 1 | 0 | 0 | 0 | 0 | 0 | 0 | 0 | 0 | 0 |
| Biomass Generator | 1 | 0 | 0 | 0 | 0 | 79 | 0 | 0 | 0 | 0 |
| Gas Co-generator-1 | 1 | 5.7 | 0 | 0 | 0 | 599 | 0 | 0 | 0 | 0 |
| Gas Co-generator-2 | 1 | 5.7 | 0 | 0 | 0 | 599 | 0 | 0 | 0 | 0 |
| Gas Co-generator-3 | 1 | 5.7 | 0 | 0 | 0 | 599 | 0 | 0 | 0 | 0 |
| Oil Co-generator -1 | 1 | 2.5 | 0 | 0 | 0 | 738 | 0 | 0 | 0 | 0 |
| Oil Co-generator -2 | 1 | 2.5 | 0 | 0 | 0 | 738 | 0 | 0 | 0 | 0 |
| Oil Co-generator -3 | 1 | 2.5 | 0 | 0 | 0 | 738 | 0 | 0 | 0 | 0 |
| Lithium Ion Battery | 0.95 | 0 | 0 | 0 | 0 | 0 | 0 | 0 | 0 | 0 |
| Lithium Ion Battery-1 | 0.95 | 0 | 0 | 0 | 0 | 0 | 0 | 0 | 0 | 0 |
| Heat Pump | 0 | 1.45 | 0 | 0 | 0 | 0 | 0 | 0 | 0 | 0 |
| Thermal Energy Storage | 0 | 1 | 0 | 0 | 0 | 0 | 0 | 0 | 0 | 0 |
| Reciprocating Internal Combustion Engine-1 | 1 | 7.06 | 0 | 0 | 0 | 670.9 | 0 | 0 | 0 | 0 |
| Reciprocating Internal Combustion Engine-2 | 1 | 4.04 | 0 | 0 | 0 | 491.7 | 0 | 0 | 0 | 0 |
| Reciprocating Internal Combustion Engine-3 | 1 | 3.02 | 0 | 0 | 0 | 448.2 | 0 | 0 | 0 | 0 |
| Gas/Combustion Turbine-1 | 1 | 6.32 | 0 | 0 | 0 | 361.5 | 0 | 0 | 0 | 0 |
| Gas/Combustion Turbine-2 | 1 | 5.14 | 0 | 0 | 0 | 302.1 | 0 | 0 | 0 | 0 |
| Gas/Combustion Turbine-3 | 1 | 5.54 | 0 | 0 | 0 | 313.4 | 0 | 0 | 0 | 0 |
| Steam Turbine-1 | 1 | 41.9 | 0 | 0 | 0 | 247.6 | 0 | 0 | 0 | 0 |
| Steam Turbine-2 | 1 | 54.5 | 0 | 0 | 0 | 247.6 | 0 | 0 | 0 | 0 |
| Steam Turbine-3 | 1 | 35.6 | 0 | 0 | 0 | 247.6 | 0 | 0 | 0 | 0 |
| Micro Turbine-1 | 1 | 5.625 | 0 | 0 | 0 | 339.7 | 0 | 0 | 0 | 0 |
| Micro Turbine-2 | 1 | 5.07 | 0 | 0 | 0 | 327.5 | 0 | 0 | 0 | 0 |
| Micro Turbine-3 | 1 | 4.93 | 0 | 0 | 0 | 369.7 | 0 | 0 | 0 | 0 |
| Fuel Cell-1 | 1 | 3.33 | 0 | 0 | 0 | 235.9 | 0 | 0 | 0 | 0 |
| Fuel Cell-2 | 1 | 4.93 | 0 | 0 | 0 | 224.5 | 0 | 0 | 0 | 0 |
| Fuel Cell-3 | 1 | 2.69 | 0 | 0 | 0 | 272.2 | 0 | 0 | 0 | 0 |
| Biogasifier-1 | 1 | 0 | 0 | 0 | 0 | 106.5 | 0 | 0 | 0 | 0 |
| Biogasifier-2 | 1 | 0 | 0 | 0 | 0 | 106.5 | 0 | 0 | 0 | 0 |
| IGCC-1 | 1 | 0 | 0 | 0 | 0 | 318.8 | 0 | 0 | 0 | 0 |
| IGCC-2 | 1 | 0 | 0 | 0 | 0 | 318.8 | 0 | 0 | 0 | 0 |
| IGCC-3 | 1 | 0 | 0 | 0 | 0 | 31.88 | 0 | 0 | 0 | 0 |
| Electrolyzer | 0 | 0 | 0 | 0 | 0 | 0 | 0 | 0 | 1 | 0 |
| Methanation Reactor | 0.002 | 0 | 0 | 1 | 0 | 0 | 0 | 0 | 0 | 0 |